\documentclass[APA,LATO1COL]{WileyNJD-v2}

\articletype{}

\historydates{}
\doiheadtext{}

\usepackage{amsmath}
\usepackage{amssymb}
\usepackage{stmaryrd}
\usepackage{bm}
\usepackage{color}
\usepackage{graphicx}

\DeclareMathOperator*{\essinf}{essinf}

\DeclareMathOperator*{\argmax}{argmax}

\newcommand{\RR}{\mathbb{R}}
\newcommand{\NN}{\mathbb{N}}

\newcommand{\EE}{\mathbb{E}}
\newcommand{\PP}{\mathbb{P}}

\newcommand{\Var}{\mathrm{Var}}
\newcommand{\Cov}{\mathrm{Cov}}

\newcommand{\sd}{\,\mathrm{d}}
\newcommand{\bw}{\mathbf{w}}

\newtheorem{example}{Example}
\newtheorem{Algorithm}{Algorithm}

\begin{document}

\title{Sampling Sup-Normalized Spectral Functions for Brown--Resnick Processes}

\author[1]{Marco Oesting}

\author[2]{Martin Schlather}

\author[2]{Claudia Schillings}

\authormark{Oesting \textsc{et al}}

\address[1]{\orgdiv{Department Mathematik}, \orgname{Universit\"at Siegen}, \orgaddress{\country{Germany}}}

\address[2]{\orgdiv{Institut f\"ur Mathematik}, \orgname{Universit\"at Mannheim}, \orgaddress{\country{Germany}}}

\corres{*Marco Oesting, Department Mathematik, Universit\"at Siegen, Walter-Flex-Str.~3, D-57072 Siegen, Germany. \email{oesting@mathematik.uni-siegen.de}}


\abstract[Summary]{Sup-normalized spectral functions form building blocks
	of max-stable and Pareto processes and therefore play an important role
	in modeling spatial extremes. For one of the most popular examples,
	the Brown--Resnick process, simulation is not straightforward. In this paper,
	we generalize two approaches for simulation via Markov Chain Monte Carlo 
	methods and rejection sampling by introducing new classes of proposal densities. In both cases, we provide an optimal choice of the proposal density with respect to sampling efficiency. The performance
	of the procedures is demonstrated in an example.}


\maketitle

\section{Introduction}

Spatial and spatio-temporal extreme value analysis aims at investigating 
extremes of quantities described by stochastic processes. 
In the classical setting, the real-valued process of interest 
$X = \{X(t), \, t \in K\}$ is sample-continuous on a compact domain 
$K \subset \RR^d$. Analysis of its extremes is often based on results of the limiting behavior of maxima of independent copies $X_i$, $i \in \NN$. Provided
that there exist continuous normalizing functions $a_n: K \to (0,\infty)$ and
$b_n: K \to \RR$ such that the process of normalized maxima 
$\{\max_{i=1}^n a_n^{-1}(t) \cdot (X_i(t)-b_n(t)), \, t \in K\}$
converges in distribution to some sample-continuous process $Z$ with 
nondegenerate margins as $n \to \infty$, the limit process $Z$ is necessarily
max-stable and we say that $X$ is in the max-domain of attraction of $Z$.

From univariate extreme value theory, it follows that the marginal distributions
of $Z$ are necessarily generalized extreme value (GEV) distributions \citep[cf.][for instance]{dhf06}. As 
max-stability is preserved under marginal transformations between different GEV 
distributions, without loss of generality, it can be assumed that $Z$ has 
standard Fr\'echet margins, i.e.\
$\PP(Z(t) \leq z) = \exp(-1/z)$, $z > 0$, for all $t \in K$. By 
\citet{dehaan-84}, any sample-continuous max-stable process with standard Fr\'echet margins can 
be represented as
\begin{align} \label{eq:gen-spec-repr}
 Z(t) =_d \max_{i \in \NN} \left\{ U_i \cdot V_i(t) \right\}, \quad t \in K,
\end{align}
where the so-called spectral processes $V_i$, $i \in \NN$, are independent 
copies of a nonnegative sample continuous stochastic process $V$ on $K$ 
satisfying $\EE\{V_i(t)\}=1$ for all $t \in K$, and 
$\sum_{i \in \NN} \delta_{U_i}$ is a Poisson point process on $(0,\infty)$ 
which is independent of the $V_i$ and has intensity measure $\Lambda$ given by
$\Lambda\{(u,\infty)\} = u^{-1}$ for all $u>0$.

Due to its complex structure, many characteristics of the max-stable process 
$Z$ in \eqref{eq:gen-spec-repr} cannot be calculated analytically, but need to
be assessed via simulations. In order to simulate $Z$ efficiently, 
\citet{OSZ17} suggest to make use of the sup-normalized spectral representation 
\begin{align} \label{eq:norm-spec-repr}
 Z(t) =_d \max_{i \in \NN} \Big\{ U_i \cdot c_\infty \cdot \frac{V^{\max}_i(t)}{\|V^{\max}_i\|_\infty} \Big\}, \quad t \in K,
\end{align}
where the $U_i$ are the same as above, the processes $V^{\max}_i$ are
independently and identically distributed, independently of the $U_i$, with
distribution $\PP(V^{\max} \in \cdot)$ given by
\begin{align} \label{eq:vmax} 
 \PP(V^{\max} \in B) = c_\infty^{-1} \cdot \int\nolimits_B \|v\|_\infty \PP(V \in \sd v),
 \quad B \in \mathcal{C}(K),
\end{align}
and $\|f\|_\infty = \sup_{t \in K} f(t)$ for every $f \in C(K)$, where $C(K)$
denotes the set of all real-valued continuous functions on $K$ equipped with the supremum norm $\|\cdot\|_\infty$ and corresponding $\sigma$-algebra $\mathcal{C}(K)$.
Here, the normalizing constant $c_\infty = \EE\{\|V\|_\infty\}$ is the so-called 
extremal coefficient of the max-stable process $Z$ over the domain $K$. In a
simulation study, \citet{OSZ17} demonstrate that simulation based on the 
sup-normalized spectral representation is competitive to other state-of-the-art
algorithms such as simulation based on extremal functions \citep{deo16} 
provided that the normalized spectral process $V^{\max}$ can be simulated efficiently.
\medskip

The law of the processes $V^{\max}$ also occurs when analyzing the extremes of 
a stochastic process $X$ in an alternative way focusing on exceedances over a
high threshold: If $X$ is in the max-domain of attraction of the max-stable 
process $Z$ in \eqref{eq:gen-spec-repr}, we have
\begin{align*}
\mathcal{L}\left( x^{-1} X(\cdot) \mid \|X\|_\infty > x\right) \longrightarrow_w 
\mathcal{L}\left( P \frac{V^{\max}(\cdot)}{\|V^{\max}\|_\infty} \right)
\end{align*}
as $x \to \infty$, where $P$ is a standard Pareto random variable and $V^{\max}$ 
is an independent process with law given in \eqref{eq:vmax}. The limit process 
$P V^{\max}(\cdot)/\|V^{\max}\|_\infty $ is called Pareto process 
\citep[cf.][]{ferreira-dehaan-14,dombry-ribatet-15}.
\medskip

Arising as sup-normalized spectral process for both max-stable and Pareto 
processes, the process $V^{\max}$ plays an important role in modelling and 
analyzing spatial extremes. As a crucial building block of spatial and spatio-temporal models, this process needs to be simulated in an efficient
way. Due to the measure transformation in \eqref{eq:vmax}, however, sampling of $V^{\max}$ is not straightforward  even in cases where the underlying spectral process $V$ can be simulated easily.

In the present paper, we focus on the simulation of $V^{\max}$ for the very popular class of log Gaussian spectral processes,
i.e.\ $V(t) = \exp(W(t))$ for some Gaussian process $W$ such that 
$\EE\{\exp(W(t))\} = 1$ for all $t \in K$. The resulting subclass of max-stable 
processes $Z$ in \eqref{eq:gen-spec-repr} comprises the only possible nontrivial 
limits of normalized maxima of rescaled Gaussian processes, the class of 
Brown--Resnick processes \citep{KSH09,kabluchko11}. In order to obtain 
Brown--Resnick processes that can be extended to stationary processes on
$\RR^d$, \cite{KSH09} consider
$ W(t) = G(t) - \Var\{G(t)\}/2$, $t \in K$,
with $G$ being a centered Gaussian process on $\RR^d$ with stationary 
increments and variogram
\begin{align*}
\textstyle \gamma(h) = \frac 1 2 \EE\left\{(G(t+h) - G(t))^2\right\}, \quad t, h \in \RR^d.
\end{align*}
It is important to note that the law of the resulting max-stable process $Z$ 
does not depend on the variance of $W$, but only on $\gamma$. Therefore, $Z$ is
called Brown--Resnick process associated to the variogram $\gamma$.

Recently, \citet{ho-dombry-17} introduced a two-step procedure to simulate the
corresponding sup-normalized process 
\begin{align*}
\frac{V^{\max}(\cdot)}{\|V^{\max}\|_\infty} = \exp\left(W^{\max}(\cdot) - \|W^{\max}\|_\infty\right)
\end{align*}
efficiently if the finite domain $K = \{t_1, \ldots, t_N\}$ is of small or moderate size: 
\begin{enumerate}
\item Sample the index $i$ of the component where the vector 
 $\mathbf{V}^{\max} = (V^{\max}(t_k))_{k=1,\ldots,N}$ assumes its maximum, i.e.\
 select one of the events $\mathbf{V}^{\max} \in S_i = \{\mathbf{s} \in (0,\infty)^N: \
 \|\mathbf{s}\|_\infty = s_i \}$, $i=1,\ldots,N$. Provided that the covariance 
 matrix $\mathbf{C}$ of the Gaussian vector $\mathbf{W} = (W(t_k))_{i=k}^N$ is 
 nonsingular, we have that this index is a.s.~unique and that the 
 probabilities of the corresponding events can be calculated in terms of the
 matrix $\mathbf{Q} \in \RR^{N \times N}$ and the vector $\mathbf{m} \in \RR^N$ given by
\begin{align*}
 \mathbf{Q} = \mathbf{C}^{-1} - \frac{\mathbf{C}^{-1} \mathbf{1}_N \mathbf{1}_N^\top \mathbf{C}^{-1}}{\mathbf{1}_N^\top \mathbf{C}^{-1} \mathbf{1}_N}
\qquad \text{and} \qquad
 \mathbf{m} = - \Big( \frac 1 2 \bm \sigma + \frac{1-\frac 1 2 \bm \sigma^\top \mathbf{C}^{-1} \mathbf{1}_N}{\mathbf{1}_N^\top \mathbf{C}^{-1} \mathbf{1}_N} \mathbf{1}_N^\top\Big) \mathbf{C}^{-1}
\end{align*}
 where $\bm \sigma = (\Var(W(t_k)))_{k=1,\ldots,N}$ is the variance vector of $\mathbf{W}$ and
 $\mathbf{1}_N = (1,\ldots,1)^\top \in \RR^N$. More precisely, by \citet{ho-dombry-17},
\begin{align*}
 \PP\left(\mathbf{V}^{\max} \in S_i\right) 
          = \frac{\det(\mathbf{Q}_{-i})^{-1/2} \exp\{\frac 1 2 \mathbf{m}_{-i}^\top \mathbf{Q}_{-i}^{-1} \mathbf{m}_{-i}\} \Phi_{N-1}(\mathbf{0}_{N-1}; \mathbf{Q}_{-i}^{-1} \mathbf{m}_{-i}, \mathbf{Q}_{-i}^{-1})}
    {\sum_{j=1}^N \det(\mathbf{Q}_{-j})^{-1/2} \exp\{\frac 1 2 \mathbf{m}_{-j}^\top \mathbf{Q}_{-j}^{-1} \mathbf{m}_{-j}\} \Phi_{N-1}(\mathbf{0}_{N-1}; \mathbf{Q}_{-j}^{-1} \mathbf{m}_{-j}, \mathbf{Q}_{-j}^{-1})},
\end{align*}
 where $\mathbf{m}_{-j}$ denotes the vector $\mathbf{m}$ after removing the $j$th 
 component, $\mathbf{Q}_{-j}$ denotes the matrix $\mathbf{Q}$ after removing the
 $j$th row and $j$th column and $\Phi_{N-1}(\mathbf{0}_{N-1}; \bm \mu, \bm
 \Sigma)$ is the distribution function of an $(N-1)$-dimensional Gaussian
 distribution with  mean vector $\bm \mu \in \RR^{N-1}$ and covariance matrix 
 $\bm \Sigma \in \RR^{(N-1) \times (N-1)}$ evaluated at 
 $\mathbf{0}_{N-1} = (0,\ldots,0) \in \RR^{N-1}$. 
\item Conditional on $\mathbf{V}^{\max} \in S_i$, we have
 $V^{\max}(t_i) / \|\mathbf{V}^{\max}\|_\infty = 1$ and the distribution of the
 vector $\mathbf{M} = (\log(V^{\max}(t_j)))_{j \neq i} - \log(\|\mathbf{V}^{\max}\|)$ is an 
 $(N-1)$-dimensional Gaussian distribution with mean vector
 $\mathbf{Q}_{-i}^{-1} \mathbf{m}_{-i}$ and covariance matrix $\mathbf{Q}_{-i}^{-1}$
 conditional on $\bm M$ being nonpositive. 
\end{enumerate} 
 However, the first step includes computationally expensive operations such as
 the evaluation of $(N-1)$-dimensional Gaussian distribution functions and the
 inversion of matrices of sizes $N \times N$ and $(N-1) \times (N-1)$. Furthermore, an efficient implementation of the second step is not straightforward. Thus, the procedure is feasible for small or moderate $N$ only.

In this paper, we will introduce alternative procedures for the simulation of $V^{\max}$, or, equivalently, $W^{\max} = \log V^{\max}$, that are
supposed to work for larger $N$, as well. To this end, we will modify a Markov
Chain Monte Carlo (MCMC) algorithm proposed by \citet{OSZ17} and a rejection
sampling approach based on ideas of \citet{dFD17}. Both procedures have 
originally been designed to sample sup-normalized spectral functions in general.
Here, we will adapt them to the specific case of Brown--Resnick processes.

\section{Simulating $W^{\max}$ via MCMC algorithms} \label{sec:MCMC}

Based on the Brown--Resnick process as our main example, we consider a 
max-stable process $Z$ with spectral process $V = e^W$ for some 
sample-continuous process $W$. Henceforth, we will always assume that the 
simulation domain $K = \{t_1,\ldots,t_N\} \subset \RR^d$ is finite and that
the corresponding spectral vector $\mathbf{W}$ possesses a density $f$ w.r.t.\ 
some measure $\mu$ on $\RR^N$.
Then, by \eqref{eq:vmax}, the transformed spectral vector $\mathbf{W}^{\max} = \log(\mathbf{V}^{\max})$, 
where the logarithm is applied componentwise, has the multivariate density
\begin{align*}
 f_{\max}(\bw) ={}& c_\infty^{-1} \max_{i=1}^N \exp(w_i) f(\bw),
  \qquad \bw \in \RR^N,
\end{align*}
which obviously has the same support as $f$, i.e.\ $\mathrm{supp}(f_{\max}) = \mathrm{supp}(f)$.

As direct sampling from the density $f_{\max}$ is rather sophisticated and the
normalizing constant $c_\infty$ is not readily available, it is quite appealing 
to choose an MCMC approach for simulation. In the present paper, we focus on Metropolis-Hastings algorithms with independence sampler \citep[cf.][for example]{tierney-94}. Denoting the strictly
positive proposal density on $\mathrm{supp}(f)$ by $f_{\rm prop}$, the algorithm
is of the following form:
\begin{Algorithm}\label{algo:simu-MCMC} MCMC APPROACH (METROPOLIS--HASTINGS)\\[1mm]
  \begin{tabular}[h]{l}\hline
 \textbf{Input:} proposal density $f_{\rm prop}$\\
 Simulate $\bw^{(0)}$ according to the density $f_{\rm prop}$.\\
 for $k=1,\ldots, n_{MCMC}$ \{\\
 \phantom{for} Sample $\bw$ from $f_{\rm prop}$ and set\\
 \phantom{for}    $
    \bw^{(k)} = \begin{cases}
        \bw         & \text{with probability }     \alpha(\bw^{(k-1)}, \bw),\\
        \bw^{(k-1)} & \text{with probability } 1 - \alpha(\bw^{(k-1)}, \bw),
      \end{cases}
                      $\\
 \phantom{for} where the acceptance probability $\alpha(\cdot,\cdot)$ is given 
  by \eqref{eq:acceptance}.\\               
\}\\
   \textbf{Output:} Markov chain $(\bw^{(1)},\ldots, \bw^{(n_{MCMC})})$.
    \\\hline
  \end{tabular}
\end{Algorithm}
\medskip

Here, the acceptance ratio $\alpha(\widetilde \bw, \bw)$ for a new proposal 
$\bw \in \mathrm{supp}(f)$ given a current state $\widetilde \bw \in 
\mathrm{supp}(f)$ is
\begin{align} \label{eq:acceptance}
 \alpha(\widetilde \bw, \bw) = \min\left\{ \frac{f_{\max}(\bw) / f_{\rm prop}(\bw)}{f_{\max}(\widetilde \bw) / f_{\rm prop}(\widetilde \bw)}, 1\right\}, 
\end{align}
using the convention that a ratio is interpreted as $0$ if both the enumerator and
the denominator are equal to $0$. This choice of $\alpha(\widetilde \bw, \bw)$  ensures reversibility of the resulting Markov chain 
$\{\bw^{(k)}\}_{k \in \NN}$ with respect to the distribution of $\mathbf{W}^{\max}$. Further, it allows for a direct transition from any state $\widetilde \bw \in \mathrm{supp}(f)$ to any other state $\bw \in \mathrm{supp}(f)$. Consequently, the chain is irreducible and aperiodic and, thus, its distribution converges to the desired stationary distribution, that is, for a.e.\ initial state 
$\bw^{(0)} \in \mathrm{supp}(f)$, we have that
\begin{align} \label{eq:conv-mcmc}
  \|P^n(\bw^{(0)}, \cdot) - \PP(\mathbf{W}^{\max} \in \cdot)\|_{\rm TV} \stackrel{n \to \infty}{\longrightarrow} 0,
\end{align}
where $P^n(\bw^{(0)}, \cdot)$ denotes the distribution of the $n$-th state of a 
Markov chain with initial state $\bw^{(0)}$ and $\|\cdot\|_{\rm TV}$ is the 
total variation norm. 
\medskip

As a general approach for the simulation of sup-normalized spectral processes 
of arbitrary max-stable processes, \citet{OSZ17} propose to use Algorithm
\ref{algo:simu-MCMC} with the density $f$ of the original spectral vector 
$\mathbf{W}$ as proposal density (Algorithm 1A) and the Metropolis-Hastings acceptance ratio in \eqref{eq:acceptance} simplifies to
\begin{align} \label{eq:acceptance-osz}
\alpha(\widetilde \bw, \bw) = \min\left\{ \frac{\max_{i=1}^N e^{w_i}}{ \max_{i=1}^N e^{\widetilde w_i}}, 1\right\}, \qquad \bw, \widetilde \bw \in \RR^N.
\end{align}
As the proposal density $f_{\rm prop} = f$ is strictly positive on
$\mathrm{supp}(f)$, convergence of the distribution of the Markov chain to 
the distribution of $\mathbf{W}^{\max}$ as in \eqref{eq:conv-mcmc} is ensured. If the support of $f$ is unbounded, however, there is no uniform
geometric rate of convergence of the chain in \eqref{eq:conv-mcmc}, as we have
\begin{align*}
\textstyle \essinf_{\bw \in \RR^N} \frac{f(\bw)}{f_{\max}(\bw)} = \essinf_{\bw \in \RR^N}
  \left(c_\infty \cdot \min_{i=1}^N e^{-w_i}\right) = 0
\end{align*}  
\citep{mengersen-tweedie-96}. In particular, this holds true for the case of
a Brown--Resnick process where $\mathbf{W}$ is a Gaussian vector.

Furthermore, due to the structure of the acceptance ratio in 
\eqref{eq:acceptance-osz}, the Markov chain may get stuck, once a state 
$\widetilde \bw$ with a large maximum 
$\max_{i=1}^N \exp(\widetilde w_i)$ is reached. This might lead to 
rather poor mixing properties of the chain. Even though independent 
realizations could still be obtained by starting new independent Markov chains
\citep[cf.][]{OSZ17}, such a behavior is undesirable having chains in high dimension $N$ with potentially long burn-in periods in mind.
\medskip 

While the algorithm in \citet{OSZ17} is designed to be applicable in a general
framework, we will use a specific transformation to construct a Markov chain 
with stronger mixing and faster convergence to the target distribution.
For many popular models such as Brown--Resnick processes, this transformation is
easily applicable. More precisely, we consider the related densities $f_i$, 
$i=1,\ldots,N$, with $f_i(\bw) = \exp(w_i) f(\bw)$. These densities are 
closely related to the distributions $P_i$ that have been studied in
\citet{deo16}.

Hence, we propose to approach the target distribution with density
 $f_{\max} = c_\infty^{-1}
\max_{i=1}^N f_i$ by Algorithm \ref{algo:simu-MCMC} using a mixture 
\begin{align} \label{eq:mixture}
f_{{\rm prop}} = \sum\nolimits_{i=1}^N p_i f_i
\end{align}
as proposal density, where the weights $p_i \geq 0$, $i=1,\ldots,N$, are such
that $\sum_{i=1}^N p_i = 1$. The corresponding acceptance probability 
in \eqref{eq:acceptance} is then
\begin{align} \label{eq:acceptance-new}
\widetilde \alpha (\widetilde \bw, \bw)
 = \min\left\{ \frac{\max_{i=1}^N e^{w_i} / \sum_{i=1}^N p_i e^{w_i}}{\max_{i=1}^N e^{\widetilde w_i} / \sum_{i=1}^N p_i e^{\widetilde w_i}}, 1\right\}.
\end{align}

With the proposal density being strictly positive on $\mathrm{supp}(f)$, we
can see that the distribution of the Markov chain again converges to its 
stationary distribution with density $f_{\max}$. As we further have
\begin{align} \label{eq:dens-ratio}
 \inf_{\bw \in \RR^N} \frac{f_{\rm prop}(\bw)}{f_{\max}(\bw)} = \inf_{\bw \in \RR^N} \frac{\sum_{i=1}^N p_i e^{w_i}}{c_\infty^{-1} \max_{i=1}^N e^{w_i}} 
 = c_\infty \cdot \min_{i=1}^N p_i > 0,
\end{align}
provided that $p_i > 0$ for $i=1,\ldots,N$,
the results found by \citet{mengersen-tweedie-96} even ensure a uniform
geometric rate of convergence for any starting value $\bw^{(0)} \in 
\mathrm{supp}(f)$ in contrast to the case where $f_{\rm prop}=f$. 
\medskip

In order to obtain a chain with good mixing properties, we choose $p_i$
such that the acceptance rate in Algorithm \ref{algo:simu-MCMC} is high
provided that the current state $\bw^{(k)}$ is (approximately) distributed
according to the stationary distribution. To this end, we 
minimize the relative deviation between $f_{\rm prop}$ and $f_{\max}$
under $f_{\max}$, i.e.\ we minimize
\begin{align*}
  D(p_1,\ldots,p_N) = \int\nolimits_{\RR^N} \bigg( \frac{f_{\rm prop}(\bw)}{f_{\max}(\bw)} - 1\bigg)^2 \, f_{\max}(\bw) \sd \bw,
\end{align*}
under the constraint $\sum_{i=1}^N p_i = 1$.
Introducing a Lagrange multiplier $\lambda \in \RR$, minimizing 
\begin{align*}
 D(p_1,\ldots,p_N) ={}& \EE\bigg\{ \bigg( \frac{\sum_{i=1}^N p_i e^{W(t_i)}}{c_\infty^{-1} e^{\max_{j=1}^N W(t_j)}} - 1\bigg)^2
       c_\infty^{-1} e^{\max_{j=1}^N W(t_j)} \bigg\} 
  {}={} c_\infty \sum_{i=1}^N \sum_{k=1}^N p_i p_k \EE\left\{e^{W(t_i) + W(t_k) - \max_{j=1}^N W(t_j)}\right\} - 1
\end{align*}
results in solving the linear system
\begin{align} \label{eq:optim-ls}
   \left( \begin{array}{cc} \bm \Sigma & \mathbf{1}_N \\ \mathbf{1}_N^\top & 0\end{array}\right) \,
   \left( \begin{array}{c} \mathbf{p}\\ \lambda \end{array}\right) = \left( \begin{array}{c} \mathbf{0}_N \\ 1 \end{array}\right),
\end{align}
where $\mathbf{p} = (p_1,\ldots,p_N)^\top$ and $\bm \Sigma = (\sigma_{ik})_{1 \leq i,k \leq N}$
with
\begin{align} \label{eq:optim-matrix-1}
 \textstyle \sigma_{ik} = \EE\big\{e^{W(t_i) + W(t_k) - \max_{j=1}^N W(t_j)}\big\}.
\end{align}
Provided that the matrix $\bm \Sigma$ is nonsingular, the solution of
\eqref{eq:optim-ls} is given by
\begin{align} \label{eq:weights-solution}
  \mathbf{p} = \frac{ \bm \Sigma^{-1} \mathbf{1}_N }{ \mathbf{1}_N^\top \bm \Sigma^{-1} \mathbf{1}_N } 
\end{align}
\citep[cf.][for instance]{cressie93}.
This solution does not necessarily satisfy the additional restriction 
$p_i \geq 0$ for all $i=1,\ldots,N$. In case that $\Sigma$ is singular or
the vector $\mathbf{p}$ has 
at least one negative entry, the full optimization problem
\begin{align} \label{eq:QP}
 \min{} & \mathbf{p}^\top \bm \Sigma \mathbf{p} \nonumber \\
 \text{s.t.} \quad \mathbf{1}_N^\top \mathbf{p} & = 1 \tag{QP} \\
                             p_i & \geq 0  \qquad \forall i =1,\ldots,N, \nonumber
\end{align}
has to be solved. Using the Karush--Kuhn--Tucker optimality conditions, the 
quadratic program \eqref{eq:QP} can be transformed into a linear program with 
additional (nonlinear) complementary slackness conditions. It can be solved 
by modified simplex methods. Alternatively, the problem \eqref{eq:QP} can be
solved by the dual method by \citet{goldfarb-idnani-1983}.
\medskip

\begin{remark}
  In order to ensure a geometric rate of convergence of the distribution of the Markov chain, we might replace the condition $p_i \geq 0$ for all $i=1,\ldots,N,$ in \eqref{eq:QP} by $p_i \geq \varepsilon$ for some given $\varepsilon > 0$. Then, a geometric rate of convergence follows from \eqref{eq:dens-ratio} as described above.
\end{remark}

In the case of Brown--Resnick processes, for simplicity, we consider the case 
that $\mathbf{W}$ possesses a full Lebesgue density. In particular, the covariance
matrix $\mathbf{C} = \left( \Cov(W(t_i), W(t_j)) \right)_{1 \leq t_i, t_j \leq N}$
of $\mathbf{W}$ is assumed to be nonsingular. Then, the target density is
\begin{align*}
 f_{\max}(\bw) {}={} c_\infty^{-1} \max_{i=1}^N \exp(w_i) f(\bw) 
               {}={} \frac{c_\infty^{-1} \max_{i=1}^N \exp(w_i)}{(2 \pi)^{N/2} \det(\mathbf{C})^{1/2}} 
    \exp\left\{ - \frac 1 2 \left(\bw + \frac{\bm \sigma} 2\right)^\top \mathbf{C}^{-1} \left(\bw + \frac{\bm \sigma}2\right)\right\}, \quad \bw \in \RR^N,
\end{align*}   
 where $\bm \sigma = (\Var(W(t_k)))_{k=1,\ldots,N}$ is again the variance vector of $\mathbf{W}$.
Now, the densities $f_i(\bw)) = e^{w_i} f(\bw)$ which form the
proposal density are just shifted Gaussian distributions:
\begin{align*}
 f_{i}(\bw) {}={} \frac {\exp(w_i)}{(2 \pi)^{\frac N 2} \det(\mathbf{C})^{\frac 1 2}} 
       \exp\left\{ - \frac 1 2 \left(\bw + \frac{\bm \sigma} 2\right)^\top \mathbf{C}^{-1} \left(\bw + \frac{\bm \sigma}2\right)\right\} 
       {}={} \frac {1}{(2 \pi)^{\frac N 2} \det(\mathbf{C})^{\frac 1 2}}
       \exp\left\{ - \frac 1 2  \left(\bw - \mathbf{C}_{\cdot i} + \frac{\bm \sigma} 2\right)^\top \mathbf{C}^{-1} \left(\bw - \mathbf{C}_{\cdot i} + \frac{\bm \sigma}2\right)\right\},
\end{align*}
cf.\ Lemma 1 in the Supplementary Material of \citet{deo16}, i.e.\
we have
\begin{align} \label{eq:shiftBR}
 \mathcal{L}(\mathbf{W}^{(i)}) = \mathcal{L}(\mathbf{W} + C_{\cdot i})
\end{align}
where the Gaussian vectors $\mathbf{W}^{(i)}$ and $\mathbf{W}$ possess densities $f_i$ 
and $f$, respectively. The calculation of the optimal weights $p_i$ is based on
the expectation in \eqref{eq:optim-matrix-1} which typically cannot be 
calculated analytically, but needs to be assessed numerically via simulations.
Such a numerical evaluation, however, is challenging as the random variable 
$\exp(W(t_i) + W(t_k) - \max_{j=1}^N W(t_j))$ is unbounded. To circumvent these
computational difficulties, we make use of the identity
\begin{align}
  \sigma_{ik} {}={} \EE\left\{e^{W(t_i) + W(t_k) - \max_{j=1}^N W(t_j)}\right\} 
  ={}& \int\nolimits_{\RR^N} \frac{e^{w_i + w_k - \max_{j=1}^N w_j}}{(2 \pi)^{\frac N 2} \det(\mathbf{C})^{\frac 1 2}} 
       \exp\left\{ - \frac 1 2 \left(\bw + \frac{\bm \sigma} 2\right)^\top \mathbf{C}^{-1} \left(\bw + \frac{\bm \sigma}2\right)\right\} \sd \bw \nonumber \\
       ={}& \int\nolimits_{\RR^N} \frac{e^{w_k - \max_{j=1}^N w_j}}{(2 \pi)^{\frac N 2} \det(\mathbf{C})^{\frac 1 2}}
       \exp\left\{ - \frac 1 2  \left(\bw - \mathbf{C}_{\cdot i} + \frac{\bm \sigma} 2\right)^\top \mathbf{C}^{-1} 
                     \left(\bw - \mathbf{C}_{\cdot i} + \frac{\bm \sigma}2\right)\right\} \sd \bw \nonumber \\
       ={}& \int\nolimits_{\RR^N} e^{w_k - \max_{j=1}^N w_j} f_i(\bw) \sd \bw
       {}={} \EE\bigg\{ \exp( W^{(i)}(t_k) - \max_{j=1}^N W^{(i)}(t_j)) \bigg\} \label{eq:sigma-calc}
\end{align}
This expression can be conveniently assessed numerically as the random variable
$\exp(W^{(i)}(t_k) - \max_{j=1}^N W^{(i)}(t_j))$ is bounded by $1$.

\begin{remark}
 Note that both \eqref{eq:shiftBR} and the final result in \eqref{eq:sigma-calc} 
 still hold true if $\mathbf{W}$ does not possess a full Lebesgue density, but exactly
 one component $W_{i^*}$ is degenerate and the reduced covariance matrix
 $(C_{ij})_{i,j \neq i^*}$ is nonsingular. This situation appears in several 
 examples such as $W$ being a fractional Brownian motion where $W(0)=0$ a.s. 
\end{remark}

In summary, we propose the procedure below to simulate the normalized 
spectral vector $\mathbf{W}^{\max}$ for the Brown--Resnick process (Algorithm 1B):
\begin{enumerate}
  \item Calculate $\mathbf{p}$ by solving the quadratic program \eqref{eq:QP}
        where the entries of the matrix $\bm \Sigma$ are given by \eqref{eq:sigma-calc}. Provided that all its components are
        nonnegative, the solution $\mathbf{p}$ has the form 
        \eqref{eq:weights-solution}.
  \item Run Algorithm \ref{algo:simu-MCMC} with proposal density 
        $f_{\rm prop} = \sum_{i=1}^N p_i f_i$ and acceptance probability
        given by \eqref{eq:acceptance-new}.
        The output of the algorithm is a Markov chain whose stationary distribution
        is the distribution of $\mathbf{W}^{\max}$. 
\end{enumerate}

\section{Exact Simulation via Rejection Sampling} \label{sec:rej-sampling}

In this section, we present an alternative procedure to generate samples
from $\mathbf{W}^{\max}$ with probability density $f_{\max}$. In contrast to Section 
\ref{sec:MCMC} where we generated a Markov chain consisting of dependent
samples with the desired distribution as stationary distribution, here, we
aim to produce independent realizations from the exact target distribution. 
To this end, we make use of a rejection sampling approach \citep[cf.][for instance]{devroye-1986} based on a proposal density $\tilde f_{\rm prop}$ 
satisfying
\begin{align} \label{eq:bound-rej-sampling}
f_{\max}(\bw) \leq ( c_\infty \cdot C)^{-1} \tilde f_{\rm prop}(\bw), \qquad 
\text{for all } \bw \in \RR^N,
\end{align}
for some $C>0$.

\begin{Algorithm}\label{algo:simu-rs} REJECTION SAMPLING APPROACH\\[1mm]
	\begin{tabular}[h]{l}\hline
		\textbf{Input:} proposal density $\tilde f_{\rm prop}$ and constant $C>0$ 
		satisfying \eqref{eq:bound-rej-sampling}\\
		repeat \{\\
		\phantom{repeat} Simulate $\bw^*$ according to the density $\tilde f_{\rm prop}$.\\
		\phantom{repeat} Generate a uniform random number $u$ in $[0,1]$.\\
		\} until $ u \cdot \tilde f_{\rm prop}(\bw^*) \leq  C \cdot c_\infty \cdot f_{\max}(\bw^*)$\\		
		\textbf{Output:} exact sample $\bw^*$ from distribution with density $f_{\max}$
		\\\hline
	\end{tabular}
\end{Algorithm}
\medskip

Thus, on average, $(c_\infty \cdot C)^{-1}$ simulations from the
proposal distribution are needed to obtain an exact sample from the target distribution. Of course, to minimize the computational burden, for a given
proposal density $\tilde f_{\rm prop}$, the constant $C$ should be chosen maximal
subject to \eqref{eq:bound-rej-sampling}, i.e.\
\begin{align*}
C = \inf_{\bw \in \RR^N} \frac{\tilde f_{\rm prop}(\bw)}{c_\infty f_{\max}(\bw)}.
\end{align*}

Recently, \citet{dFD17} followed a similar idea and suggested to base the simulation of a general sup-normalized spectral process $V(\cdot)^{\max} / \|V^{\max}\|_\infty$ on the relation
\begin{align} \label{eq:def-dav}
 \PP\left( \frac{V^{\max}}{\|V^{\max}\|_\infty} \in \mathrm{d}v \right)
   = \frac{\|\tilde V\|_\infty}{\EE \|\tilde V\|_\infty} 
   \PP\left( \frac{\tilde V}{\|\tilde V\|_\infty} \in \mathrm{d}v \right) 
\end{align}   
where $\tilde V$ is a spectral process normalized with respect to another
homogeneous functional $r$ instead of the supremum norm, i.e.\ 
$r(\tilde V) = 1$ a.s.
 If $\|\tilde V\|_\infty$ is a.s.\ bounded from
above by some constant, from the relation \eqref{eq:def-dav}, we obtain
an inequality of the same type as \eqref{eq:bound-rej-sampling} for the 
densities of $V(\cdot)^{\max} / \|V^{\max}\|_\infty$ and $\tilde V(\cdot) / \|\tilde V\|_\infty$ instead of $f_{\max}$ and $f_{\rm prop}$, respectively.
Thus, samples of $\tilde V(\cdot) / \|\tilde V\|_\infty$ can be used as 
proposals for an exact rejection sampling procedure. For instance, the sum-normalized spectral vector $\tilde{\mathbf{V}}$, i.e.\ the vector which is normalized w.r.t.\ the functional $r(f) =  \|f\|_1 = \sum_{k=1}^N |f(t_k)|$,
can be chosen as it is easy to simulate in many cases \citep[cf.][]{deo16}
and satisfies $\|\tilde{\mathbf{V}}\|_\infty \leq 1$ almost surely.
\medskip

For a Brown--Resnick process, it is well-known that the sum-normalized process
$\tilde{\mathbf{V}}$ has the same distribution as
$\exp(\mathbf{W}^{\rm prop}) / \|\exp(\mathbf{W}^{\rm prop})\|_1$
where $\mathbf{W}^{\rm prop}$ has the density $f_{\rm prop}$ from \eqref{eq:mixture}
in Section \ref{sec:MCMC} with equal weights $p_1 = \ldots = p_N = 1/N$
\citep[see also][]{dieker-mikosch-15}. Thus, 
in this case, the procedure proposed by \citet{dFD17} with $r(f) = \|f\|_1$
is equivalent to performing rejection sampling for $\mathbf{W}^{\max}$ with 
$\tilde f_{\rm prop} = \sum_{i=1}^N \frac 1 N f_i$ as proposal distribution
(Algorithm 2A). From Equation \eqref{eq:dens-ratio}, 
it follows that rejection sampling can also be performed with $\tilde f_{\rm prop} = \sum_{i=1}^N p_i f_i$ and arbitrary positive 
weights $p_1,\ldots,p_N$ summing up to $1$, since we have 
\eqref{eq:bound-rej-sampling} with $C= \min_{i=1}^N p_i$.
Thus, accepting a proposal $\bw^*$ in the rejection sampling procedure with
probability
\begin{align*} 
 \min_{i=1}^N p_i \cdot \frac{c_{\infty} f_{\max}(\bw^*)}{f_{\rm prop}(\bw^*)} 
= \frac{\min_{i=1}^N p_i \cdot \max_{i=1}^N e^{w_i^*}}
       {\sum_{i=1}^N p_i e^{w_i^*}},
\end{align*}       
we will obtain a sample of independent realizations from the exact target
distribution $f_{\max}$. The rejection rate, however, is pretty high. In order
to obtain one realization of $\mathbf{W}^{\max}$, on average 
$\left(c_\infty \cdot \min_{i=1}^N p_i\right)^{-1}$
simulations from $f^*_{\rm prop}$ are required. It can be easily seen that the
computational costs are indeed minimal for the choice $p_1=\ldots=p_N = 1/N$,
i.e.\ the choice in the approach based on the sum-normalized representation.
In this case, one realization of $\mathbf{W}^{\max}$ on average requires to sample
$c_\infty^{-1} N$ times from $f_{\rm prop}$. Therefore, this approach becomes
rather inefficient if we have a large number $N$ of points on a dense grid.
\medskip

In order to reduce the large computational costs of rejection sampling which 
are mainly due to the fact that $\min_{i=1}^N f_i(\bw) / f_{\max}(\bw)$ 
gets small as $\|\bw\| \to \infty$ , we replace each density $f_i$ by the 
modified multivariate Gaussian density $g_{i,\varepsilon}$ whose variance is
increased by the factor $(1-\varepsilon)^{-1} \geq 1$ for some 
$\varepsilon \in [0,1)$:
\begin{align*}
  g_{i,\varepsilon}(\bw) 
     {}={}& \frac{(1-\varepsilon)^{N/2}}{(2 \pi)^{N/2} \det(\mathbf{C})^{1/2}} \exp\left( - \frac 1 2 (1-\varepsilon) 
               \cdot \Big(\bw - \mathbf{C}_{\cdot i} + \frac{\bm \sigma} 2\Big)^\top \mathbf{C}^{-1} \Big(\bw - \mathbf{C}_{\cdot i} + \frac{\bm \sigma} 2\Big) \right) \\
     ={}& \frac{(1-\varepsilon)^{N/2}}{(2 \pi)^{N/2} \det(\mathbf{C})^{1/2}} \exp\left( (1-\varepsilon) w_i\right)
               \cdot \exp\left( - \frac 1 2 (1-\varepsilon) \cdot \Big(\bw + \frac{\bm \sigma} 2\Big)^\top \mathbf{C}^{-1} \Big(\bw + \frac{\bm \sigma} 2\Big) \right).
\end{align*}
Analogously to $f_{\rm prop}$ for the MCMC approach in Section \ref{sec:MCMC},
we propose a mixture 
\begin{align*}
\tilde f_{\rm prop} = \sum\nolimits_{i=1}^N p_i g_{i,\varepsilon} 
\end{align*}
with $p_i \geq 0$ and $\sum_{i=1}^N p_i = 1$ as proposal density
for the rejection sampling algorithm. A proposal $\bw^*$ is then
accepted with probability
\begin{align} \label{eq:acc-g}
  C(\mathbf{p}, \varepsilon) \cdot \frac{c_\infty \cdot f_{\max}(\bw^*)}
        {\sum_{i=1}^N p_i g_{i,\varepsilon}(\bw^*)}
 {}={}& C(\mathbf{p}, \varepsilon) \cdot 
 \frac{\exp\Big(-\frac {\varepsilon} 2 \Big(\bw + \frac{\bm \sigma} 2\Big)^\top \mathbf{C}^{-1} \Big(\bw + \frac{\bm \sigma}{2}\Big) \Big)}
    { (1-\varepsilon)^{\frac N 2} \cdot \sum_{i=1}^N p_i \exp\left( (1-\varepsilon) w_i - \max_{j=1}^N w^*_j \right)},
\end{align}
where
\begin{align} \label{eq:def-c}
  C(\mathbf{p}, \varepsilon) 
 ={}& \inf_{\bw \in \RR^N} \frac{\sum_{i=1}^N p_i g_{i,\varepsilon}(\bw)}{c_\infty f_{\max}(\bw)} {}={} \inf_{\bw \in \RR^N} \min_{j=1}^N \frac{\sum_{i=1}^N p_i g_{i,\varepsilon}(\bw)}{f_j(\bw)} \nonumber \\
 ={}& \inf_{\bw \in \RR^N} (1-\varepsilon)^{\frac N 2} \sum_{i=1}^N p_i
    \exp\left( (1-\varepsilon) w_i - \max_{j=1}^N w_j + \frac {\varepsilon} 2 \Big(\bw + \frac{\bm \sigma} 2\Big)^\top \mathbf{C}^{-1} \Big(\bw + \frac{\bm \sigma}{2}\Big) \right).
\end{align}

Thus, to summarize, for appropriately chosen
$\varepsilon > 0$ and $\mathbf{p} \geq 0$ such that $\|\mathbf{p}\|_1=1$,
we propose to run  Algorithm \ref{algo:simu-rs} with proposal density
$ \tilde f_{\rm prop} = \sum\nolimits_{i=1}^N p_i g_{i,\varepsilon} $ 
and $C = C(\mathbf{p},\varepsilon)$ according to \eqref{eq:def-c}.

\begin{remark}
	To further reduce the computational costs in the simulation, we might even
	choose a more flexible approach. For instance, instead of using a mixture of a
	finite number of functions $g_{1,\varepsilon},\ldots,g_{N,\varepsilon}$, one
	could consider arbitrary mixtures
	\begin{align*}
	\tilde f_{\rm prop}(\bw) = \int\nolimits_{\RR^d} g_{t,\varepsilon}(\bw) \nu(\mathrm{d}t), \qquad \bw \in \RR^N,
	\end{align*}
	where $g_{t,\varepsilon}(\bw) = \int_{\RR} g_{N+1,\varepsilon}(\bw, w_{N+1}) \mathrm{d} w_{N+1}$
	on the enlarged domain $\{t_1,\ldots,t_N, t\}$ and $\nu$ is a probability 
	measure on $\RR^d$. Furthermore, depending on $t \in \RR^d$, different values
	for $\varepsilon=\varepsilon(t) \in [0,1)$ might be chosen. However, due to 
	the complexity of the optimization problems involved, we restrict ourselves to
	the situation above where $\nu$ is a probability measure on 
	$K = \{t_1,\ldots,t_N\}$ and $\varepsilon$ is constant in space.
\end{remark}

Using the procedure described above, on average, $(c_\infty \cdot C(\mathbf{p}, \varepsilon))^{-1}$ simulations from the proposal distribution are needed
to obtain one exact sample from the target distribution, i.e.\ the computational
complexity of the algorithm depends on the choices of $\mathbf{p}$ and 
$\varepsilon$. The remainder of this section will be devoted to this question.

\paragraph{Choice of $\mathbf{p}$ and $\varepsilon$} 
For a given $\varepsilon \geq 0$, the computational costs of the algorithm can be minimized by choosing $\mathbf{p} = \mathbf{p}^*(\varepsilon)$ such that the constant
$C(\mathbf{p}, \varepsilon)$ given in \eqref{eq:def-c} is maximal, i.e.\ by choosing
$\mathbf{p}$ as the solution of the nonlinear optimization problem
\begin{align}\label{eq:N}
	\max_{\mathbf{p} \in \RR^N} & \quad C(\mathbf{p}, \varepsilon) \nonumber \\
	\text{s.t.} \quad \|\mathbf{p}\|_1 & =     1 \tag{NP}\\
				p_i   & \geq  0 \qquad \forall i=1,\ldots,N\nonumber
\end{align}
Optimizing further w.r.t.\ $\varepsilon \in [0,1)$, we obtain 
the optimal choice $(\mathbf{p}, \varepsilon) = (\mathbf{p}^*(\varepsilon^*),
\varepsilon^*)$ where $\varepsilon^* = \argmax_{\varepsilon \in [0,1)} C(\mathbf{p}^*(\varepsilon), \varepsilon)$. 

As the above optimization problem includes optimization steps w.r.t.\
$\bw \in \RR^N$, $\mathbf{p} \in \{\mathbf{x} \in [0,1]^N: \ \|\mathbf{x}\|_1 = 1\}$ and
$\varepsilon \in [0,1]$, none of which can be solved analytically, the solution
is quite involved. In order to reduce the computational burden, we simplify the
problem by maximizing an analytically simpler lower bound. To this end, we
decompose the convex combination $\sum_{i=1}^N p_i g_{i,\varepsilon} / f_j$
into sums over disjoint subsets of the form $I = \{i_1,\ldots,i_m\} \subset \{1,\ldots,N\}$.
For a convex combination of $(g_{i_k,\varepsilon})_{k=1,\ldots,m}$ with weight
vector $\bm \lambda = (\lambda_k)_{k=1,\ldots,m} \in [0,1]^m$, we obtain the 
lower bound
\begin{align*}
   \inf_{\bw \in \RR^N} \frac{\sum_{k=1}^m \lambda_k g_{i_k,\varepsilon}(\bw)}{f_j(\bw)}
 ={}&  \inf_{\bw \in \RR^N} (1-\varepsilon)^{N/2} \cdot \left( \sum\nolimits_{k=1}^m \lambda_k e^{(1-\varepsilon) w_{i_k} - w_j} \right) 
         \cdot \exp\left( \frac{\varepsilon} 2 \left(\bw + \frac{\bm \sigma} 2\right)^\top \mathbf{C}^{-1} \left(\bw + \frac{\bm \sigma}{2}\right) \right) \displaybreak[0]\\
 \geq{}& \inf_{\bw \in \RR^N} (1-\varepsilon)^{N/2} \cdot \exp\left( (1-\varepsilon) \sum\nolimits_{k=1}^m \lambda_k w_{i_k} - w_j \right) 
         \cdot \exp\left( \frac{\varepsilon} 2 \left(\bw + \frac{\bm \sigma} 2\right)^\top \mathbf{C}^{-1} \left(\bw + \frac{\bm \sigma}{2}\right) \right) {}=:{} c_I^{(j)}(\varepsilon, \bm \lambda),   
\end{align*}
where we made use of the convexity of the exponential function. Setting 
$\kappa_I^{(j)}(\varepsilon,\bm \lambda) = (1-\varepsilon) \sum_{k=1}^m \lambda_k \mathbf{C}_{\cdot i_k} - \mathbf{C}_{\cdot j}$, 
this bound can be calculated explicitly:
\begin{align} \label{eq:cIj}
 c_I^{(j)}(\varepsilon, \bm \lambda)  ={}& \inf_{\bw \in \RR^N} (1-\varepsilon)^{N/2} \cdot
     \exp\left\{ \frac{\varepsilon}{2} \left(\bw + \frac 1 {\varepsilon} \kappa_I^{(j)}(\varepsilon,\bm \lambda) + \frac{\bm \sigma}{2}\right)^\top
                            \mathbf{C}^{-1} \left(\bw + \frac 1 {\varepsilon} \kappa_I^{(j)}(\varepsilon,\bm \lambda) + \frac{\bm \sigma}{2}\right) \right\} \nonumber \\
    & \hspace{2.5cm} \cdot \exp\left( - (\kappa_I^{(j)}(\varepsilon,\bm \lambda))^\top \mathbf{C}^{-1} \frac{\bm \sigma} 2 
                            - \frac 1 {2\varepsilon} (\kappa_I^{(j)}(\varepsilon,\bm \lambda))^\top \mathbf{C}^{-1} \kappa_I^{(j)}(\varepsilon,\bm \lambda) \right) \nonumber \displaybreak[0] \\
 ={}& (1-\varepsilon)^{N/2} \cdot \exp\left( - \frac{1-\varepsilon} 2 \sum\nolimits_{k=1}^m \lambda_k C_{i_k i_k} + \frac 1 2 C_{jj} 
                                             - \frac {(1-\varepsilon)^2}{2\varepsilon} \sum\nolimits_{k=1}^m \lambda_k^2 C_{i_k i_k} - \frac 1 {2\varepsilon} C_{jj} \right) \nonumber \\
    & \hspace{1.6cm} \cdot \exp\left( - \frac{(1-\varepsilon)^2}{2\varepsilon} \sum\nolimits_{k=1}^m \sum\nolimits_{l \neq k} \lambda_k \lambda_l C_{i_k i_l} 
                                    + \frac{ 1-\varepsilon   }{ \varepsilon} \sum\nolimits_{k=1}^m \lambda_k C_{i_k j}\right) \nonumber \\
 ={}& (1-\varepsilon)^{N/2} \exp\left( - \frac{1-\varepsilon}{\varepsilon} \sum\nolimits_{k=1}^m \lambda_k \gamma(t_{i_k} - t_j) 
                                       + \frac{(1-\varepsilon)^2}{2\varepsilon} \sum\nolimits_{k=1}^m \sum\nolimits_{l=1}^m \lambda_k \lambda_l \gamma(t_{i_k} - t_{i_l}) \right).
\end{align}

Hence,
\begin{align*}
 \sum\nolimits_{i \in I} p_i g_{i,\varepsilon}(\bw)
  ={}& \|\mathbf{p}_I\|_1 \cdot \sum\nolimits_{i \in I} \frac{p_i}{\|\mathbf{p}_I\|_1} g_{i,\varepsilon}(\bw) \geq{} \|\mathbf{p}_I\|_1 \cdot c_{I}^{(j)}\left(\varepsilon, \frac{\mathbf{p}_I}{\|\mathbf{p}_I\|_1}\right) \cdot f_j(\bw), 
    \qquad \bw \in \RR^N,
\end{align*}
where $\mathbf{p}_I = (p_i)_{i \in I}$ for every subset $I \subset \{1,\ldots,N\}$.
Now, for each $j \in \{1,\ldots,N\}$, let $J^{(j)}$ be a partition of 
$\{1,\ldots,N\}$, so that
 \begin{align} \label{eq:sum-bound}
   c_\infty \cdot C(\mathbf{p}, \varepsilon) 
   {}={}& \inf_{\bw \in \RR^N} \frac{g_{\rm prop}(\bw)}{f_{\max}(\bw)} 
   {}={} \inf_{\bw \in \RR^N} \min_{j=1}^N \frac{\sum_{I \in J^{(j)}} \sum_{i \in I} p_i g_{i,\varepsilon}(\bw)}{c_\infty^{-1} f_j(\bw)} \nonumber \\
   {}\geq{}& \inf_{\bw \in \RR^N} \min_{j=1}^N 
     \frac{\sum_{I \in J^{(j)}} \left(\sum\nolimits_{i \in I} p_i\right) \cdot c_{I}^{(j)}\left(\varepsilon, \frac{\mathbf{p}_I}{\|\mathbf{p}_I\|_1}\right) \cdot f_j(\bw)}
          {c_\infty^{-1} f_j(\bw)} 
   {}={} c_\infty \cdot \min_{j=1}^N \sum\nolimits_{I \in J^{(j)}} \left(\sum\nolimits_{i \in I} p_i\right) \cdot c_{I}^{(j)}\left(\varepsilon, \frac{\mathbf{p}_I}{\|\mathbf{p}_I\|_1}\right).
 \end{align}
Thus, the RHS of \eqref{eq:sum-bound} provides an explicit lower bound for the
average acceptance probability for any choice of the $J^{(j)}$.
 
\begin{remark}
Assume that, for some $j \in \{1,\ldots,N\}$, there is some index set 
$I = \{i_1,\ldots,i_m\} \subset \{1,\ldots,N\}$ such that 
$\gamma(t_{i} - t_j) = \Gamma$ for all $i \in I$.
Then, Equation \eqref{eq:cIj} provides the bound
\begin{align} \label{eq:bound-group}
       \sum_{k=1}^m p_{i_k} g_{i_k,\varepsilon}(\bw) \geq{}& \|\mathbf{p}_I\|_1 c_{I}^{(j)} \left(\varepsilon, \frac{\mathbf{p}_I}{\|\mathbf{p}_I\|_1}\right) f_j(\bw)
    ={} \|\mathbf{p}_I\|_1 \exp\left( - \frac{1-\varepsilon}{\varepsilon} \Gamma 
        +  \frac{(1-\varepsilon)^2}{2\varepsilon} \sum\nolimits_{k=1}^m \sum\nolimits_{l=1}^m \frac{p_{i_k} p_{i_l}}{\|\mathbf{p}_I\|_1^2} \gamma(t_{i_k} - t_{i_l}) \right) f_j(\bw)
\end{align}
for all $\bw \in \RR^N$. Alternatively, for the same index set $I$, we could
bound each summand separately, i.e.\
\begin{align} \label{eq:bound-single}
   \sum\nolimits_{k=1}^m p_{i_k} g_{i_k,\varepsilon}(\bw) {}\geq{} \sum\nolimits_{k=1}^m p_{i_k} c_{\{i_k\}}^{(j)}(\varepsilon,1) f_j(\bw) 
   {}={} \|\mathbf{p}_I\|_1\exp\left( - \frac{1-\varepsilon}{\varepsilon} \Gamma \right) f_j(\bw).
\end{align}
Note that, for all $\bw \in \RR^N$, the RHS of \eqref{eq:bound-single} is less
than the RHS of \eqref{eq:bound-group}, i.e.\ the lower bound is less sharp. 
Therefore, we prefer pooling locations with the same distance to $t_j$ rather 
than considering them separately in order to have the bound in 
\eqref{eq:sum-bound} as sharp as possible.
\end{remark}

In view of \eqref{eq:sum-bound}, instead of considering the exact value
$C(\mathbf{p}, \varepsilon)$ which is needed to calculate the minimal rejection
rate, but cannot be given explicitly, we might maximize the 
function
\begin{align*}
 C_{\mathrm{groups}}(\mathbf{p}, \varepsilon) 
 = \min_{j=1}^N \sum_{I \in J^{(j)}} \sum_{i \in I} p_i c_I^{(j)}\left(\varepsilon, \frac{\mathbf{p}_I}{\|\mathbf{p}_I\|}\right)
\end{align*}
for fixed partitions $J^{(1)},\ldots,J^{(N)}$. Due to the complex dependence of
$c_I^{(j)}$ on $\mathbf{p}$, the resulting optimization problem is nonlinear in $\mathbf{p}$
even for fixed $\varepsilon$. To circumvent this difficulty, for each $I \in J^{(j)}$,
we fix $|I|$-dimensional weight vectors $\bm \lambda(I)$ and consider the 
function
\begin{align*}
 C_{\mathrm{groups}}^{\mathrm{fix}}(\mathbf{p}, \varepsilon; \bm \lambda) ={}& \min_{j=1}^N \sum_{I \in J^{(j)}} \sum_{i \in I} p_i c_I^{(j)}\left(\varepsilon, \bm \lambda(I)\right)
 ={}  \min_{j=1}^N \sum\nolimits_{i=1}^N p_i c_{ij}(\varepsilon; \bm \lambda) = \min \mathbf{p}^\top \mathbf{c}(\varepsilon; \bm \lambda)
\end{align*}
with $\mathbf{c}(\varepsilon; \bm \lambda) = \{c_{ij}(\varepsilon; \bm \lambda)\}_{1 \leq i,j \leq N}$
and $c_{ij}(\varepsilon; \bm \lambda) = c_I^{(j)}(\varepsilon; \bm \lambda(I))$ 
for the unique set $I \in J^{(j)}$ such that $i \in I$. 
\medskip

Analogously to the solution above, we first maximize 
$C_{\mathrm{groups}}^{\mathrm{fix}}(\cdot, \varepsilon; \bm \lambda)$ for fixed 
$\varepsilon \in [0,1)$ and $\bm \lambda$, i.e.\ we consider the optimization 
problem
\begin{align}\label{eq:LP1}
	\max_{\mathbf{p} \in \RR^N} \min_{j=1}^N \ & c_{\cdot j}(\varepsilon; \bm \lambda)^\top \mathbf{p} \nonumber \\
	\text{s.t.} \quad \|\mathbf{p}\|_1&=1 \tag{LP1}\\
				p_i   & \geq  0 \qquad \forall i=1,\ldots,N\nonumber\,.
\end{align}
To convert the linear program to standard form, we introduce an additional
variable $z \in \RR$, unconstrained in sign, leading to the equivalent program
\begin{align} \label{eq:LP2}
	\max_{\mathbf{p} \in \RR^N, z \in \RR} {}& z  \nonumber\\
	\text{s.t.} 
	\quad z & \leq c_{\cdot j}(\varepsilon; \bm \lambda)^\top \mathbf{p} 
	\qquad \forall j=1,\ldots,N \tag{LP2}\\
	\mathbf{1}_{N}^\top \mathbf{p}&=1 \nonumber \\
    p_i  & \geq	 0 \qquad \forall i=1,\ldots,N\nonumber\,.
\end{align}
The standard form of \eqref{eq:LP2} is then given by
\begin{align}\label{eq:LP2S}
  \max_{\substack{p \in \RR^N\\ s,z^+,z^- \in \RR}}& z^+-z^- \nonumber  \\
   \text{s.t.}\quad z^+-z^-+s & = c_{\cdot j}(\varepsilon; \bm \lambda)^\top \mathbf{p} \qquad \forall j=1,\ldots,N \tag{LP2S}\\
        \mathbf{1}_{N}^\top \mathbf{p} & = 1 \nonumber \\
        p_1,\ldots,p_N, s, z^+,z^- & \geq 0 \nonumber\,.
\end{align}
Such a linear program in standard form can be solved by standard techniques such 
as the simplex algorithm. Compared to the optimization of $C_{\mathrm{groups}}(\cdot, \varepsilon)$,
the complementary one-dimensional problem of maximizing
$C_{\mathrm{groups}}(\mathbf{p}, \cdot)$ for fixed $\mathbf{p}$ can be solved rather
easily.
\medskip

To summarize, starting from some $\varepsilon^*>0$ and
 $\mathbf{p}^* = N^{-1} \mathbf{1}_N$, we propose to apply the following two steps repeatedly (Algorithm 2B):
\begin{enumerate}
 \item Define 
 \begin{align*}
 \bm \lambda_I = \frac{\mathbf{p}_I^*}{\|\mathbf{p}_I^*\|_1}, \quad I \in J^{(1)} \cup \ldots \cup J^{(N)}
 \end{align*}
    and set 
    \begin{align*}   
      \mathbf{p}^* = \argmax_{\mathbf{p}} \, C_{\mathrm{groups}}^{\mathrm{fix}}(\mathbf{p}, \varepsilon^*; \bm \lambda),
    \end{align*}
    i.e.\ the solution of the optimization problem \eqref{eq:LP1} (or 
    \eqref{eq:LP2} or \eqref{eq:LP2S}, equivalently).
 \item Set $\varepsilon^* = \argmax_{\varepsilon} C_{\mathrm{groups}}(\mathbf{p}^*, \varepsilon)$.
\end{enumerate}

Even though $C_{\mathrm{groups}}(\mathbf{p}, \varepsilon)$ might be significantly
smaller than $C(\mathbf{p}, \varepsilon)$, in some cases, this bound is already 
sufficient to improve the results for $\varepsilon = 0$ that have been
discussed in the beginning of this section, where we have already seen that 
the corresponding optimal weight vector equals $\mathbf{p}^* = N^{-1} \mathbf{1}_N$ and 
that $C(\mathbf{p}^*,0) = 1/N$. We show an example to illustrate that this choice 
is not necessarily optimal, i.e.\ there is some $\varepsilon > 0$ and a vector 
$\mathbf{p}$ of weights such that $C(\mathbf{p}, \varepsilon) \geq
C_{\mathrm{groups}}(\mathbf{p}, \varepsilon) > 1/N$.

\begin{example}[Fractional Brownian Motion] \label{ex:rej}
 Let $\mathbf{x}_1, \ldots, \mathbf{x}_N$ be $N$ equidistant locations in $[0,1]$
 and $\gamma(h) = |h|^\alpha$ for some $\alpha > 1$. Choose
 $\mathbf{p} = N^{-1} \mathbf{1}_N$, $\varepsilon = \sqrt{2} N^{-1}$
 and set $J^{(1)} = \ldots = J^{(N)} =\{\{1\},\ldots,\{N\}\}$.
 Then, for every $x_i$ there are at least $\lfloor N^{-1/\alpha} \cdot N \rfloor$ 
 locations $\mathbf{x}_j$ such that $\gamma(\mathbf{x}_i - \mathbf{x}_j) \leq 1/N$.
 Thus, we obtain for large $N$ that
 \begin{align*}
   C(\mathbf{p}, \varepsilon) \geq 
   C_{\mathrm{groups}}(\mathbf{p}, \varepsilon) \geq{} \lfloor N^{-1/\alpha} \cdot N \rfloor
   \cdot \frac 1 N \cdot
   \left( 1 - \frac{\sqrt{2}} N \right)^{N/2} 
   \exp\left( - \frac{N}{\sqrt{2}} \cdot \frac 1 N\right)
   \sim{}  N^{-1/\alpha} \exp(-\sqrt{2})
 \end{align*}
 which is eventually larger than $1/N$ as $\alpha > 1$.
\end{example}

\section{Illustration}

Finally, we illustrate the performance of Algorithm \ref{algo:simu-MCMC} and 
the rejection sampling algorithm in an example. Taking up Example \ref{ex:rej} in
higher dimension, we consider the case that $Z$ is a Brown--Resnick process
associated to the variogram
\begin{align*}
\gamma(h) = \left\| \frac{h}{5} \right\|^{1.5}, \quad h \in \RR^2,
\end{align*}
on the grid $K = \{0, 0.2,\ldots,5\} \times \{0, 0.2,\ldots,5\}$ 
($N=676$ points).
We run four different algorithms:
\begin{enumerate}
 \item[1A.] Algorithm \ref{algo:simu-MCMC} with proposal density 
      $f_{\rm prop} = f$ as proposed by \citet{OSZ17}
 \item[1B.] Algorithm \ref{algo:simu-MCMC} with proposal density 
       $f_{\rm prop} = \sum_{i=1}^N p_i f_i$ where $\mathbf{p}$ is given as the
       solution of \eqref{eq:QP} (cf.\ Section \ref{sec:MCMC}).
 \item[2A.] Algorithm \ref{algo:simu-rs} with proposal density $f_{\rm prop} = 
       \frac 1 N \sum_{i=1}^N f_i$ and $C=1/N$ (equivalent to the procedure proposed in \citet{dFD17} based on sum-normalized spectral functions)
 \item[2B.] Algorithm \ref{algo:simu-rs}  with proposal density $\tilde f_{\rm prop} = \sum_{i=1}^N p_i^*(\varepsilon^*) g_{i,\varepsilon^*}$ 
 and $C = C(\mathbf{p}^*, \varepsilon^*)$ in \eqref{eq:def-c} where 
       $\mathbf{p}^*(\varepsilon^*) \in [0,1]^N$ and $\varepsilon^* \in 
       (0,1)$ are obtained as described in Section \ref{sec:rej-sampling}
\end{enumerate}

Even though the laws of the Brown--Resnick process $Z$ and the normalized
spectral process $V^{\max}$ do not depend on the variance, but only on the 
variogram of the underlying Gaussian process $G$, the choice of the Gaussian
process may affect the performance of the algorithms. Here, we choose the
Gaussian process $G$ whose law is uniquely defined via the construction
\begin{align*}
G(t) = G_0(t) - \frac 1 4 \left( G_0\left(\begin{pmatrix} 0 \\ 0 \end{pmatrix}\right) +  G_0\left(\begin{pmatrix} 5 \\ 0 \end{pmatrix}\right) 
                                 +  G_0\left(\begin{pmatrix} 0 \\ 5 \end{pmatrix}\right) +  G_0\left(\begin{pmatrix} 5 \\ 5 \end{pmatrix}\right) \right), \quad t \in K,
\end{align*}                                 
where $G_0$ is an arbitrary centered Gaussian process with variogram $\gamma$.
\citet{oesting-strokorb-2017} show that this process has a smaller maximal
variance and is thus preferable in the context of simulation.

\begin{figure}
 \centering \includegraphics[width=7.2cm]{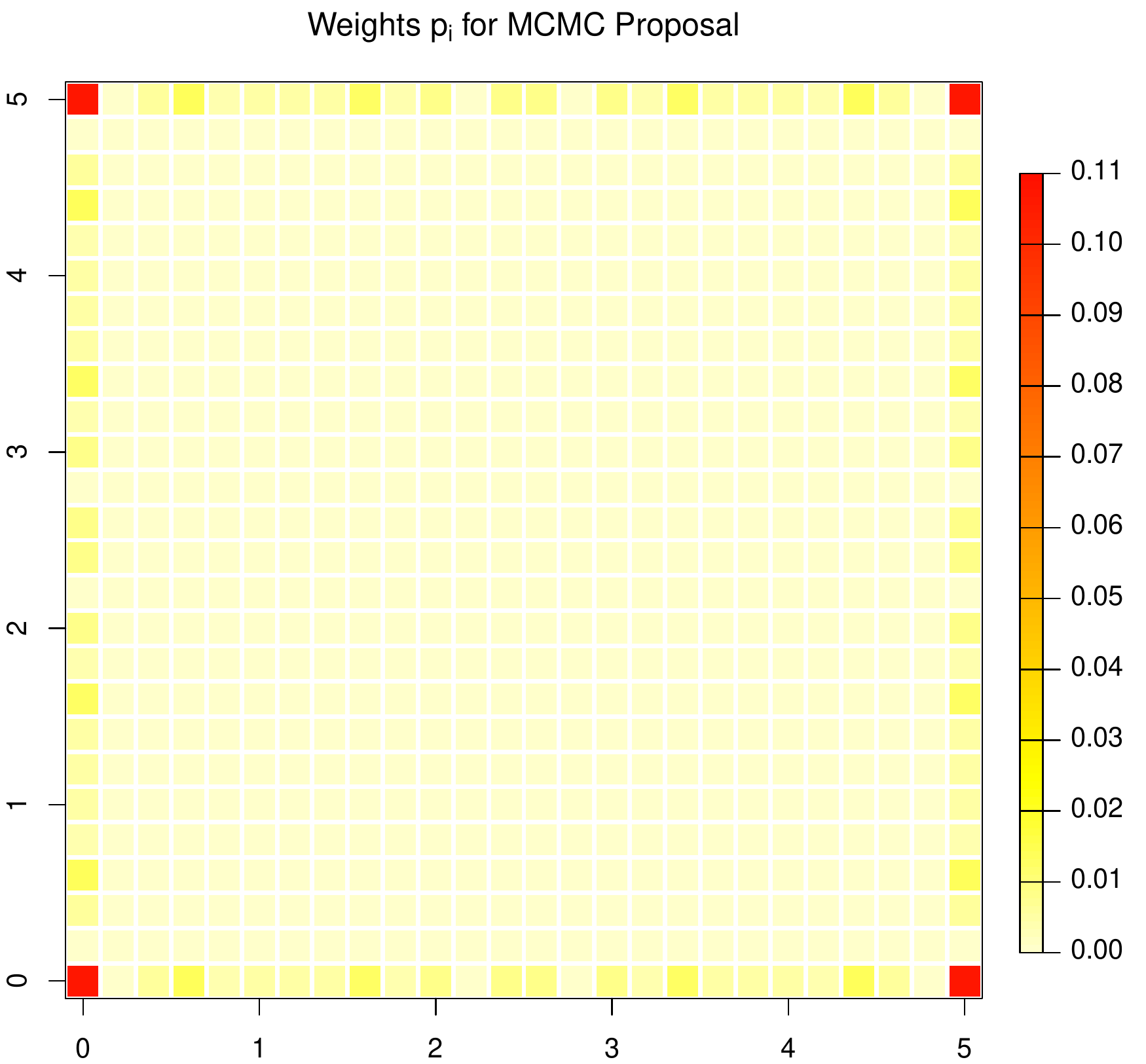}\includegraphics[width=7.2cm]{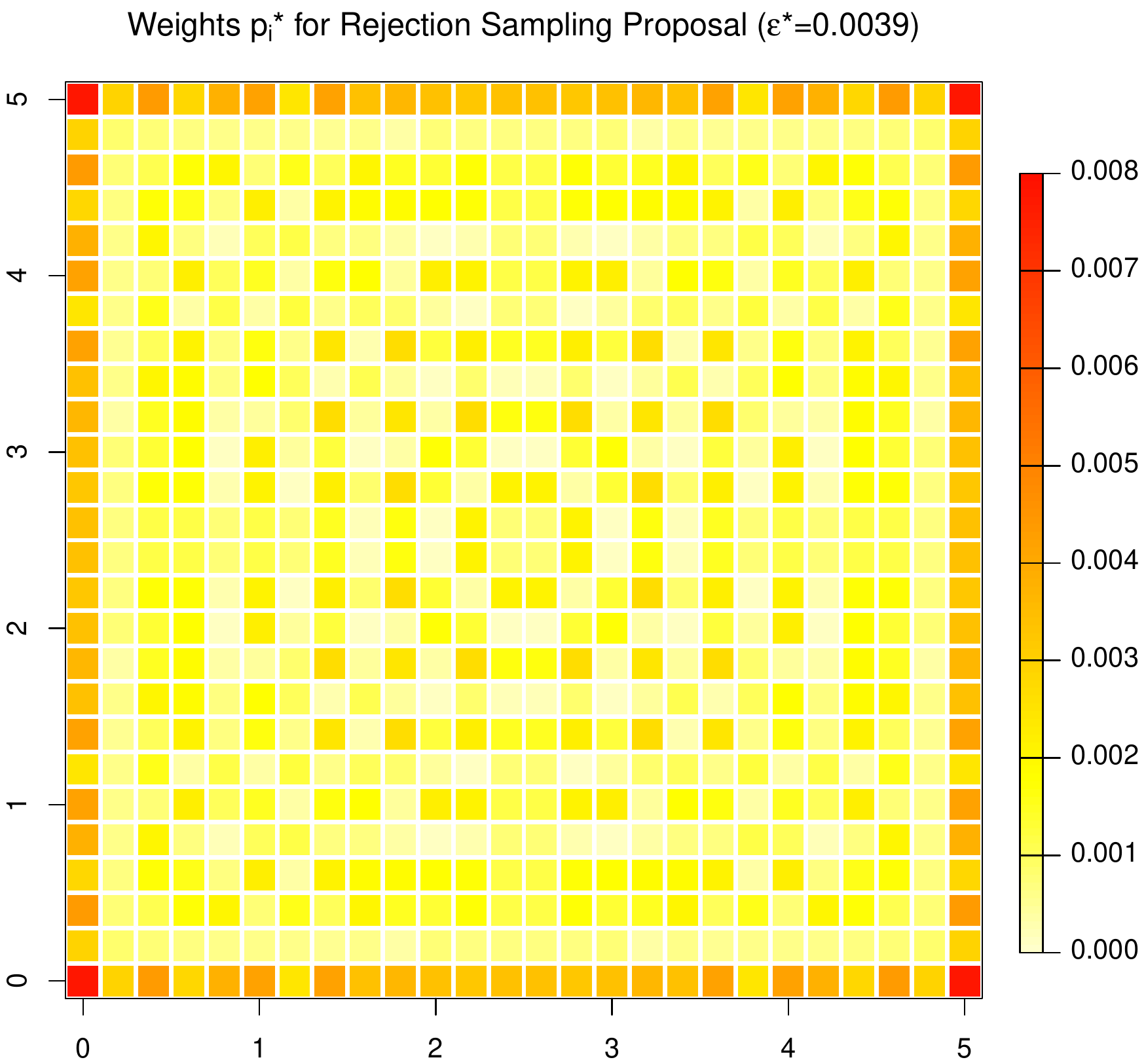}
 
 \caption{Vectors $\mathbf{p}$ (left) and $\mathbf{p}^*$ (right) of optimal weights used
   in Algorithms 1B and 2B, respectively, for a Brown--Resnick process associated 
   to the variogram $\gamma(h) = \|h/5\|^{1.5}$ on the set 
   $K = \{0, 0.2,\ldots, 5\}^2$.} \label{fig:weights}
\end{figure}

We first calculate the optimal weights $\mathbf{p} = (p_1,\ldots,p_{676})^\top$
as a solution of \eqref{eq:QP} (used in Algorithm 1B) as well as
 the optimal weights  $\mathbf{p}^*(\varepsilon^*)$ as a solution of \eqref{eq:LP1} and the optimal variance modification $\varepsilon^*$ (used in Algorithm 2B). The results for $\mathbf{p}$ and $\mathbf{p}^*(\varepsilon^*)$ are displayed in Figure \ref{fig:weights}. 
 It can be seen that, in both 
cases, the weights are not spatially constant, but are larger on the boundary
of the convex hull $\mathrm{conv}(K) = [0,5] \times [0,5]$ with the maximum in
the corners of the square. This observation is well in line with the fact that
these points have the largest contribution to 
$\max_{t \in K} \exp(G(t) - \Var(G(t))/2)$ since the variance of $G$ attains 
its maximum there \citep[see also][]{oesting-strokorb-2017}. 

\begin{figure}
 \centering \includegraphics[width=14.4cm]{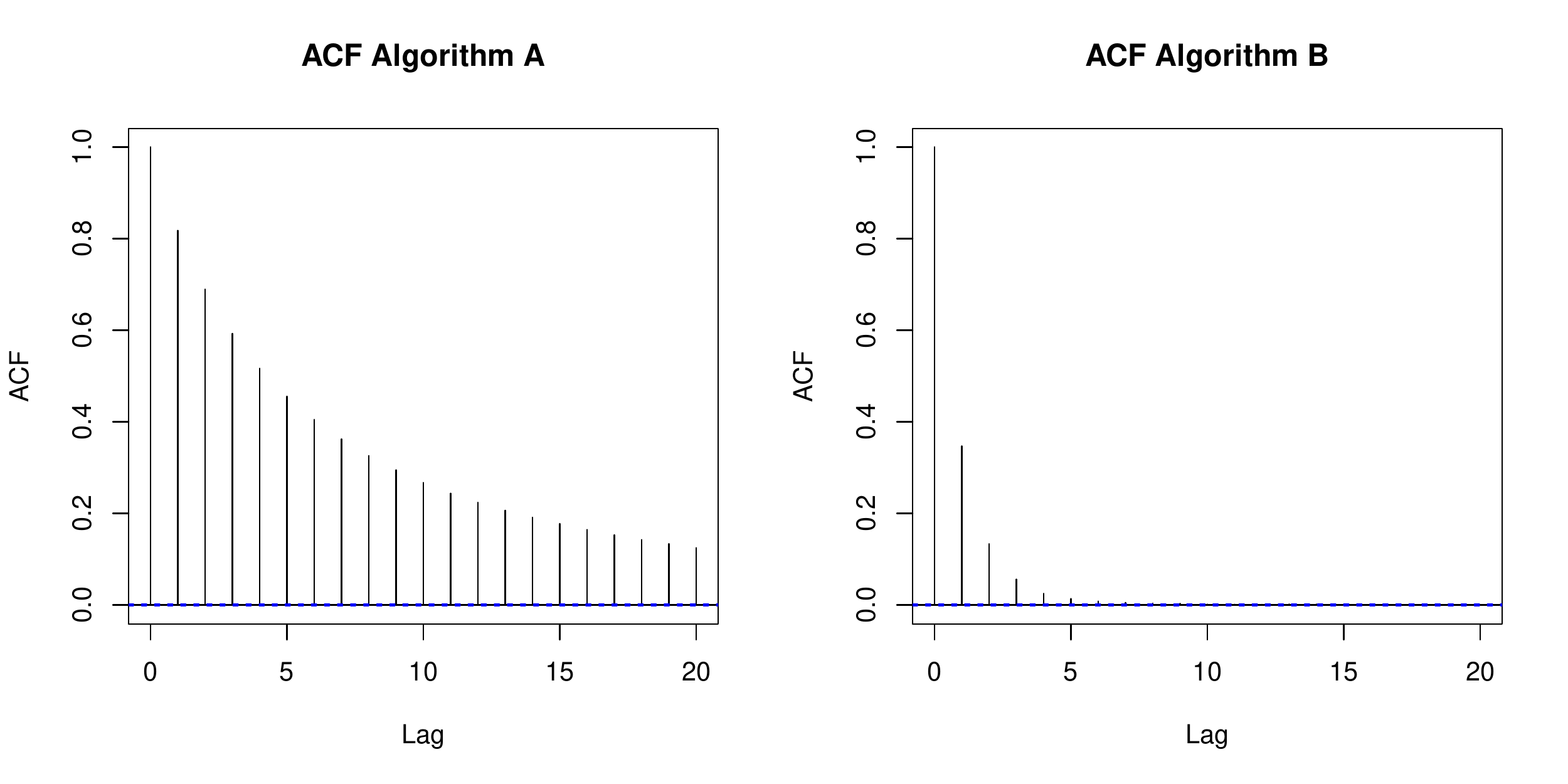}
 \caption{Empirical autocorrelation functions of the time series
     $\{\|\exp(W^{(k)}_{1})\|_\infty\}_{k=1,\ldots,1\,000\,000}$ (left)
     and $\{\|\exp(W^{(k)}_{2})\|_\infty\}_{k=1,\ldots,1\,000\,000}$ (right)
     obtained via Algorithms 1A and 1B, respectively.
     Here, by the optimal choice of $\mathbf{p}$, the autocorrelation is clearly
     reduced.
 } \label{fig:acf}
\end{figure}

We run Algorithm \ref{algo:simu-MCMC} with both proposal densities as specified
above (Algorithms 1A and 1B, respectively) to obtain two Markov chains 
$\{W^{(k)}_{1}\}_{k=1,\ldots,1\,000\,000}$ and $\{W^{(k)}_{2}\}_{k=1,\ldots,1\,000\,000}$
of length $n_{MCMC} = 1\,000\,000$. It can be seen that the empirical acceptance
rate of the second chain ($0.855$) is remarkably higher than the one of the first
chain ($0.656$) which already indicates stronger mixing. This impression is
confirmed by analyzing the empirical autocorrelation functions of the time series 
$\{\|\exp(W^{(k)}_{1})\|_\infty\}_{k=1,\ldots,1\,000\,000}$ and  
$\{\|\exp(W^{(k)}_{2})\|_\infty\}_{k=1,\ldots,1\,000\,000}$ which are shown in 
Figure \ref{fig:acf}. Here, the empirical autocorrelation in the second chain
is drastically reduced in comparison with the first chain, indicating that two
states of the chain can be regarded as nearly uncorrelated after roughly five steps.

The rejection sampling algorithm (Algorithm 2A and Algorithm 2B) automatically generates independent realizations from the multivariate
target density $f_{\max}$. Therefore, we will compare them with 
respect to their computational complexity. To this end, we run them to generate
a sample of size $100\,000$ and count the average number of simulations of 
Gaussian vectors from the proposal density to generate one realization of 
$\mathbf{W}^{\max}$. In case of Algorithm 2A, this number is $203.1$ which is close
to the theoretical expression $c_\infty^{-1} \cdot N$. For Algorithm 2B, the
number is improved by a factor of approximately $4.4$, leading to an average
number of $45.9$ Gaussian vectors to be simulated to obtain one realization 
from the target distribution. This improvement is well in line with the
corresponding value $C(\mathbf{p}^*, \varepsilon^*) \approx 0.0065 \approx 4.4 \cdot
N^{-1}$.
\medskip

As the example illustrates, the two modifications we suggested may lead to 
significant improvements of MCMC and rejection algorithms that have been 
proposed so far. Here, only the modified rejection sampling algorithm ensures
independence of exact samples from the target distribution. However, as the 
example indicates, the MCMC algorithm might be particularly attractive in 
practice as a thinned chain results in nearly independent samples even if
the thinning rate is rather small. 
Note that we also tried other examples such as a Brownian sheet ($\alpha=1$). However, we found that significant improvements in the rejection sampling procedure become apparent only for $\alpha > 1$, see also Example \ref{ex:rej}.

\section*{Acknowledgements}

The authors are grateful to Kirstin Strokorb, Dimitri Schwab and Jonas Brehmer for numerous valuable comments.  M.~Schlather has been financially supported by Volkswagen Stiftung within the project `Mesoscale Weather Extremes -- Theory, Spatial Modeling and Prediction (WEX-MOP)'. 

\bibliography{lit.bib}

\begin{thebibliography}{}

\bibitem [\protect \citeauthoryear {%
Cressie%
}{%
Cressie%
}{%
{\protect \APACyear {1993}}%
}]{%
cressie93}
\APACinsertmetastar {%
cressie93}%
\begin{APACrefauthors}%
Cressie, N\BPBI A.%
\end{APACrefauthors}%
\unskip\
\newblock
\APACrefYear{1993}.
\newblock
\APACrefbtitle {Statistics for {S}patial {D}ata} {Statistics for {S}patial
  {D}ata}.
\newblock
\APACaddressPublisher{New York}{John Wiley \& Sons}.
\PrintBackRefs{\CurrentBib}

\bibitem [\protect \citeauthoryear {%
de Fondeville%
\ \BBA {} Davison%
}{%
de Fondeville%
\ \BBA {} Davison%
}{%
{\protect \APACyear {2018}}%
}]{%
dFD17}
\APACinsertmetastar {%
dFD17}%
\begin{APACrefauthors}%
de Fondeville, R.%
\BCBT {}\ \BBA {} Davison, A\BPBI C.%
\end{APACrefauthors}%
\unskip\
\newblock
\APACrefYearMonthDay{2018}{}{}.
\newblock
{\BBOQ}\APACrefatitle {High-dimensional peaks-over-threshold inference}
  {High-dimensional peaks-over-threshold inference}.{\BBCQ}
\newblock
\APACjournalVolNumPages{Biometrika}{105}{3}{575-592}.
\PrintBackRefs{\CurrentBib}

\bibitem [\protect \citeauthoryear {%
de Haan%
}{%
de Haan%
}{%
{\protect \APACyear {1984}}%
}]{%
dehaan-84}
\APACinsertmetastar {%
dehaan-84}%
\begin{APACrefauthors}%
de Haan, L.%
\end{APACrefauthors}%
\unskip\
\newblock
\APACrefYearMonthDay{1984}{}{}.
\newblock
{\BBOQ}\APACrefatitle {A spectral representation for max-stable processes} {A
  spectral representation for max-stable processes}.{\BBCQ}
\newblock
\APACjournalVolNumPages{Ann. Probab.}{12}{4}{1194--1204}.
\PrintBackRefs{\CurrentBib}

\bibitem [\protect \citeauthoryear {%
de Haan%
\ \BBA {} Ferreira%
}{%
de Haan%
\ \BBA {} Ferreira%
}{%
{\protect \APACyear {2006}}%
}]{%
dhf06}
\APACinsertmetastar {%
dhf06}%
\begin{APACrefauthors}%
de Haan, L.%
\BCBT {}\ \BBA {} Ferreira, A.%
\end{APACrefauthors}%
\unskip\
\newblock
\APACrefYear{2006}.
\newblock
\APACrefbtitle {Extreme {V}alue {T}heory: {A}n {I}ntroduction} {Extreme {V}alue
  {T}heory: {A}n {I}ntroduction}.
\newblock
\APACaddressPublisher{New York}{Springer}.
\PrintBackRefs{\CurrentBib}

\bibitem [\protect \citeauthoryear {%
Devroye%
}{%
Devroye%
}{%
{\protect \APACyear {1986}}%
}]{%
devroye-1986}
\APACinsertmetastar {%
devroye-1986}%
\begin{APACrefauthors}%
Devroye, L.%
\end{APACrefauthors}%
\unskip\
\newblock
\APACrefYear{1986}.
\newblock
\APACrefbtitle {{N}on-{U}niform {R}andom {V}ariate {G}eneration}
  {{N}on-{U}niform {R}andom {V}ariate {G}eneration}.
\newblock
\APACaddressPublisher{New York}{Springer-Verlag}.
\PrintBackRefs{\CurrentBib}

\bibitem [\protect \citeauthoryear {%
Dieker%
\ \BBA {} Mikosch%
}{%
Dieker%
\ \BBA {} Mikosch%
}{%
{\protect \APACyear {2015}}%
}]{%
dieker-mikosch-15}
\APACinsertmetastar {%
dieker-mikosch-15}%
\begin{APACrefauthors}%
Dieker, A\BPBI B.%
\BCBT {}\ \BBA {} Mikosch, T.%
\end{APACrefauthors}%
\unskip\
\newblock
\APACrefYearMonthDay{2015}{}{}.
\newblock
{\BBOQ}\APACrefatitle {Exact simulation of {B}rown-{R}esnick random fields at a
  finite number of locations} {Exact simulation of {B}rown-{R}esnick random
  fields at a finite number of locations}.{\BBCQ}
\newblock
\APACjournalVolNumPages{Extremes}{18}{2}{301--314}.
\PrintBackRefs{\CurrentBib}

\bibitem [\protect \citeauthoryear {%
Dombry%
, Engelke%
\BCBL {}\ \BBA {} Oesting%
}{%
Dombry%
\ \protect \BOthers {.}}{%
{\protect \APACyear {2016}}%
}]{%
deo16}
\APACinsertmetastar {%
deo16}%
\begin{APACrefauthors}%
Dombry, C.%
, Engelke, S.%
\BCBL {}\ \BBA {} Oesting, M.%
\end{APACrefauthors}%
\unskip\
\newblock
\APACrefYearMonthDay{2016}{}{}.
\newblock
{\BBOQ}\APACrefatitle {Exact Simulation of Max-Stable Processes} {Exact
  simulation of max-stable processes}.{\BBCQ}
\newblock
\APACjournalVolNumPages{Biometrika}{103}{2}{303--317}.
\PrintBackRefs{\CurrentBib}

\bibitem [\protect \citeauthoryear {%
Dombry%
\ \BBA {} Ribatet%
}{%
Dombry%
\ \BBA {} Ribatet%
}{%
{\protect \APACyear {2015}}%
}]{%
dombry-ribatet-15}
\APACinsertmetastar {%
dombry-ribatet-15}%
\begin{APACrefauthors}%
Dombry, C.%
\BCBT {}\ \BBA {} Ribatet, M.%
\end{APACrefauthors}%
\unskip\
\newblock
\APACrefYearMonthDay{2015}{}{}.
\newblock
{\BBOQ}\APACrefatitle {Functional regular variations, {P}areto processes and
  peaks over threshold} {Functional regular variations, {P}areto processes and
  peaks over threshold}.{\BBCQ}
\newblock
\APACjournalVolNumPages{Stat. Interface}{8}{1}{9--17}.
\PrintBackRefs{\CurrentBib}

\bibitem [\protect \citeauthoryear {%
Ferreira%
\ \BBA {} de Haan%
}{%
Ferreira%
\ \BBA {} de Haan%
}{%
{\protect \APACyear {2014}}%
}]{%
ferreira-dehaan-14}
\APACinsertmetastar {%
ferreira-dehaan-14}%
\begin{APACrefauthors}%
Ferreira, A.%
\BCBT {}\ \BBA {} de Haan, L.%
\end{APACrefauthors}%
\unskip\
\newblock
\APACrefYearMonthDay{2014}{}{}.
\newblock
{\BBOQ}\APACrefatitle {The generalized {P}areto process; with a view towards
  application and simulation} {The generalized {P}areto process; with a view
  towards application and simulation}.{\BBCQ}
\newblock
\APACjournalVolNumPages{Bernoulli}{20}{4}{1717--1737}.
\PrintBackRefs{\CurrentBib}

\bibitem [\protect \citeauthoryear {%
Goldfarb%
\ \BBA {} Idnani%
}{%
Goldfarb%
\ \BBA {} Idnani%
}{%
{\protect \APACyear {1983}}%
}]{%
goldfarb-idnani-1983}
\APACinsertmetastar {%
goldfarb-idnani-1983}%
\begin{APACrefauthors}%
Goldfarb, D.%
\BCBT {}\ \BBA {} Idnani, A.%
\end{APACrefauthors}%
\unskip\
\newblock
\APACrefYearMonthDay{1983}{}{}.
\newblock
{\BBOQ}\APACrefatitle {A numerically stable dual method for solving strictly
  convex quadratic programs} {A numerically stable dual method for solving
  strictly convex quadratic programs}.{\BBCQ}
\newblock
\APACjournalVolNumPages{Math. Program.}{27}{1}{1--33}.
\PrintBackRefs{\CurrentBib}

\bibitem [\protect \citeauthoryear {%
Ho%
\ \BBA {} Dombry%
}{%
Ho%
\ \BBA {} Dombry%
}{%
{\protect \APACyear {2017}}%
}]{%
ho-dombry-17}
\APACinsertmetastar {%
ho-dombry-17}%
\begin{APACrefauthors}%
Ho, Z\BPBI W\BPBI O.%
\BCBT {}\ \BBA {} Dombry, C.%
\end{APACrefauthors}%
\unskip\
\newblock
\APACrefYearMonthDay{2017}{}{}.
\newblock
\APACrefbtitle {Simple models for multivariate regular variations and the
  {H}\"usler--{R}eiss {P}areto distribution.} {Simple models for multivariate
  regular variations and the {H}\"usler--{R}eiss {P}areto distribution.}
\newblock
\APACrefnote{arXiv preprint arXiv:1712.09225}
\PrintBackRefs{\CurrentBib}

\bibitem [\protect \citeauthoryear {%
Kabluchko%
}{%
Kabluchko%
}{%
{\protect \APACyear {2011}}%
}]{%
kabluchko11}
\APACinsertmetastar {%
kabluchko11}%
\begin{APACrefauthors}%
Kabluchko, Z.%
\end{APACrefauthors}%
\unskip\
\newblock
\APACrefYearMonthDay{2011}{}{}.
\newblock
{\BBOQ}\APACrefatitle {Extremes of independent {G}aussian processes} {Extremes
  of independent {G}aussian processes}.{\BBCQ}
\newblock
\APACjournalVolNumPages{Extremes}{14}{3}{285--310}.
\PrintBackRefs{\CurrentBib}

\bibitem [\protect \citeauthoryear {%
Kabluchko%
, Schlather%
\BCBL {}\ \BBA {} de Haan%
}{%
Kabluchko%
\ \protect \BOthers {.}}{%
{\protect \APACyear {2009}}%
}]{%
KSH09}
\APACinsertmetastar {%
KSH09}%
\begin{APACrefauthors}%
Kabluchko, Z.%
, Schlather, M.%
\BCBL {}\ \BBA {} de Haan, L.%
\end{APACrefauthors}%
\unskip\
\newblock
\APACrefYearMonthDay{2009}{}{}.
\newblock
{\BBOQ}\APACrefatitle {Stationary max-stable fields associated to negative
  definite functions} {Stationary max-stable fields associated to negative
  definite functions}.{\BBCQ}
\newblock
\APACjournalVolNumPages{Ann. Probab.}{37}{5}{2042--2065}.
\PrintBackRefs{\CurrentBib}

\bibitem [\protect \citeauthoryear {%
Mengersen%
\ \BBA {} Tweedie%
}{%
Mengersen%
\ \BBA {} Tweedie%
}{%
{\protect \APACyear {1996}}%
}]{%
mengersen-tweedie-96}
\APACinsertmetastar {%
mengersen-tweedie-96}%
\begin{APACrefauthors}%
Mengersen, K\BPBI L.%
\BCBT {}\ \BBA {} Tweedie, R\BPBI L.%
\end{APACrefauthors}%
\unskip\
\newblock
\APACrefYearMonthDay{1996}{}{}.
\newblock
{\BBOQ}\APACrefatitle {Rates of convergence of the {H}astings and {M}etropolis
  algorithms} {Rates of convergence of the {H}astings and {M}etropolis
  algorithms}.{\BBCQ}
\newblock
\APACjournalVolNumPages{Ann. Stat.}{24}{1}{101--121}.
\PrintBackRefs{\CurrentBib}

\bibitem [\protect \citeauthoryear {%
Oesting%
, Schlather%
\BCBL {}\ \BBA {} Zhou%
}{%
Oesting%
\ \protect \BOthers {.}}{%
{\protect \APACyear {2018}}%
}]{%
OSZ17}
\APACinsertmetastar {%
OSZ17}%
\begin{APACrefauthors}%
Oesting, M.%
, Schlather, M.%
\BCBL {}\ \BBA {} Zhou, C.%
\end{APACrefauthors}%
\unskip\
\newblock
\APACrefYearMonthDay{2018}{}{}.
\newblock
{\BBOQ}\APACrefatitle {Exact and Fast Simulation of Max-Stable Processes on a
  Compact Set Using the Normalized Spectral Representation} {Exact and fast
  simulation of max-stable processes on a compact set using the normalized
  spectral representation}.{\BBCQ}
\newblock
\APACjournalVolNumPages{Bernoulli}{24}{2}{1497--1530}.
\PrintBackRefs{\CurrentBib}

\bibitem [\protect \citeauthoryear {%
Oesting%
\ \BBA {} Strokorb%
}{%
Oesting%
\ \BBA {} Strokorb%
}{%
{\protect \APACyear {2018}}%
}]{%
oesting-strokorb-2017}
\APACinsertmetastar {%
oesting-strokorb-2017}%
\begin{APACrefauthors}%
Oesting, M.%
\BCBT {}\ \BBA {} Strokorb, K.%
\end{APACrefauthors}%
\unskip\
\newblock
\APACrefYearMonthDay{2018}{}{}.
\newblock
{\BBOQ}\APACrefatitle {Efficient simulation of {B}rown--{R}esnick processes
  based on variance reduction of {G}aussian processes} {Efficient simulation of
  {B}rown--{R}esnick processes based on variance reduction of {G}aussian
  processes}.{\BBCQ}
\newblock
\APACjournalVolNumPages{Adv. Appl. Probab.}{50}{4}{1155--1175}.
\PrintBackRefs{\CurrentBib}

\bibitem [\protect \citeauthoryear {%
Tierney%
}{%
Tierney%
}{%
{\protect \APACyear {1994}}%
}]{%
tierney-94}
\APACinsertmetastar {%
tierney-94}%
\begin{APACrefauthors}%
Tierney, L.%
\end{APACrefauthors}%
\unskip\
\newblock
\APACrefYearMonthDay{1994}{}{}.
\newblock
{\BBOQ}\APACrefatitle {Markov chains for exploring posterior distributions}
  {Markov chains for exploring posterior distributions}.{\BBCQ}
\newblock
\APACjournalVolNumPages{Ann. Stat.}{22}{4}{1701--1728}.
\PrintBackRefs{\CurrentBib}

\end{thebibliography}

\end{document}